\documentclass{amsart}
\usepackage{amsfonts}

\newtheorem{theorem}{Theorem}
\theoremstyle{plain}

\newtheorem{proposition}{Proposition}
\newtheorem{remark}{Remark}

\numberwithin{equation}{section}
\input{tcilatex}

\begin{document}
\title[Schwarz, Triangle and Bessel Inequalities]{New Reverses of Schwarz,
Triangle and Bessel Inequalities in Inner Product Spaces}
\author{S.S. Dragomir}
\address{School of Computer Science and Mathematics\\
Victoria University of Technology\\
PO Box 14428, MCMC 8001\\
Victoria, Australia.}
\email{sever.dragomir@vu.edu.au}
\urladdr{http://rgmia.vu.edu.au/SSDragomirWeb.html}
\date{3 September, 2003.}
\subjclass[2000]{ 46C05, 26D15.}
\keywords{Schwarz's inequality, Triangle inequality, Bessel's inequality, Gr%
\"{u}ss type inequalities, Integral inequalities.}

\begin{abstract}
New reverses of the Schwarz, triangle and Bessel inequalities in inner
product spaces are pointed out. These results complement the recent ones
obtained by the author in the earlier paper \cite{SSDa}. Further, they are
employed to establish new Gr\"{u}ss type inequalities. Finally, some natural
integral inequalities are stated as well.
\end{abstract}

\maketitle

\section{Introduction}

Let $\left( H;\left\langle \cdot ,\cdot \right\rangle \right) $ be an inner
product space over the real or complex number field $\mathbb{K}$.

In earlier paper \cite{SSDa}, we have obtained the following simple reverse
of Schwarz's inequality%
\begin{align}
0& \leq \left\Vert x\right\Vert ^{2}\left\Vert a\right\Vert ^{2}-\left\vert
\left\langle x,a\right\rangle \right\vert ^{2}  \label{1.1} \\
& \leq \left\Vert x\right\Vert ^{2}\left\Vert a\right\Vert ^{2}-\left[ \func{%
Re}\left\langle x,a\right\rangle \right] ^{2}\leq r^{2}\left\Vert
x\right\Vert ^{2},  \notag
\end{align}%
provided 
\begin{equation}
\left\Vert x-a\right\Vert \leq r<\left\Vert a\right\Vert ,  \label{1.2}
\end{equation}%
where $a,x\in H$ and $r>0.$ The constant $c=1$ in front of $r^{2}$ is best
possible in the sense that it cannot be replaced by a smaller one.

This result has then been employed to prove (see \cite{SSDa}) that%
\begin{align}
\left\Vert x\right\Vert ^{2}\left\Vert y\right\Vert ^{2}& \leq \frac{1}{4}%
\cdot \frac{\left\{ \func{Re}\left[ \left( \bar{\Gamma}+\bar{\gamma}\right)
\left\langle x,y\right\rangle \right] \right\} ^{2}}{\func{Re}\left( \Gamma 
\bar{\gamma}\right) }  \label{1.3} \\
& \leq \frac{1}{4}\cdot \frac{\left\vert \Gamma +\gamma \right\vert ^{2}}{%
\func{Re}\left( \Gamma \bar{\gamma}\right) }\left\vert \left\langle
x,y\right\rangle \right\vert ^{2},  \notag
\end{align}%
provided, for $x,y\in H$ and $\gamma ,\Gamma \in \mathbb{K}$ with $\func{Re}%
\left( \Gamma \bar{\gamma}\right) >0,$ either%
\begin{equation}
\func{Re}\left\langle \Gamma y-x,x-\gamma y\right\rangle \geq 0,  \label{1.4}
\end{equation}%
or, equivalently,%
\begin{equation}
\left\Vert x-\frac{\gamma +\Gamma }{2}y\right\Vert \leq \frac{1}{2}%
\left\vert \Gamma -\gamma \right\vert \left\Vert y\right\Vert ,  \label{1.5}
\end{equation}%
holds. In both inequalities $\left( \ref{1.3}\right) $, $\frac{1}{4}$ is the
best possible constant.

The inequality (\ref{1.3}) implies the following additive version of reverse
Schwarz's inequality%
\begin{equation}
0\leq \left\Vert x\right\Vert ^{2}\left\Vert y\right\Vert ^{2}-\left\vert
\left\langle x,y\right\rangle \right\vert ^{2}\leq \frac{1}{4}\cdot \frac{%
\left\vert \Gamma -\gamma \right\vert ^{2}}{\func{Re}\left( \Gamma \bar{%
\gamma}\right) }\left\vert \left\langle x,y\right\rangle \right\vert ^{2}.
\label{1.6}
\end{equation}%
Here the constant $\frac{1}{4}$ is also the best.

If the condition (\ref{1.2}) is satisfied, one may deduce the following
reverse of the triangle inequality \cite{SSDa}%
\begin{equation}
0\leq \left\Vert x\right\Vert +\left\Vert a\right\Vert -\left\Vert
x+a\right\Vert \leq \sqrt{2}r\sqrt{\frac{\func{Re}\left\langle
x,a\right\rangle }{\sqrt{\left\Vert a\right\Vert ^{2}-r^{2}}\left( \sqrt{%
\left\Vert a\right\Vert ^{2}-r^{2}}+\left\Vert a\right\Vert \right) }.}
\label{1.7}
\end{equation}

If $M>m>0,$ $x,y\in H$ and either (\ref{1.4}) or, equivalently, (\ref{1.5})
holds for $M,m$ instead of $\Gamma ,\gamma ,$ then the following simpler
reverse of the triangle inequality may be stated as well%
\begin{equation}
0\leq \left\Vert x\right\Vert +\left\Vert y\right\Vert -\left\Vert
x+y\right\Vert \leq \frac{\sqrt{M}-\sqrt{m}}{\sqrt[4]{Mm}}\sqrt{\func{Re}%
\left\langle x,y\right\rangle }.  \label{1.8}
\end{equation}

Moving now onto Gr\"{u}ss type inequalities, we note that if $x,y,e\in H,$
with $\left\Vert e\right\Vert =1$ and $r_{1},r_{2}\in \left( 0,1\right) $
are such that%
\begin{equation}
\left\Vert x-e\right\Vert \leq r_{1},\ \ \ \ \ \ \ \ \ \ \ \left\Vert
y-e\right\Vert \leq r_{2},  \label{1.9}
\end{equation}%
then one has the inequality \cite{SSDa}%
\begin{equation}
\left\vert \left\langle x,y\right\rangle -\left\langle x,e\right\rangle
\left\langle e,y\right\rangle \right\vert \leq r_{1}r_{2}\left\Vert
x\right\Vert \left\Vert y\right\Vert .  \label{1.10}
\end{equation}%
The inequality (\ref{1.10}) is sharp in the sense that the constant $c=1$ in
front of $r_{1}r_{2}$ cannot be replaced by a smaller constant.

If we assumed that, for $x,y,e\in H$ with $\left\Vert e\right\Vert =1$ and $%
\gamma ,\Gamma \in \mathbb{K}$ with $\func{Re}\left( \Gamma \bar{\gamma}%
\right) ,$ $\func{Re}\left( \Phi \bar{\phi}\right) >0,$ either the condition%
\begin{equation}
\func{Re}\left\langle \Gamma e-x,x-\gamma e\right\rangle ,\ \func{Re}%
\left\langle \Phi e-y,y-\gamma \phi \right\rangle \geq 0  \label{1.11}
\end{equation}%
or, equivalently,%
\begin{equation}
\left\Vert x-\frac{\gamma +\Gamma }{2}e\right\Vert \leq \frac{1}{2}%
\left\vert \Gamma -\gamma \right\vert ,\ \ \ \ \left\Vert y-\frac{\phi +\Phi 
}{2}e\right\Vert \leq \frac{1}{2}\left\vert \Phi -\phi \right\vert ,
\label{1.12}
\end{equation}%
holds, then we have the inequality%
\begin{equation}
\left\vert \left\langle x,y\right\rangle -\left\langle x,e\right\rangle
\left\langle e,y\right\rangle \right\vert \leq \frac{1}{4}\cdot \frac{%
\left\vert \Gamma -\gamma \right\vert \left\vert \Phi -\phi \right\vert }{%
\sqrt{\func{Re}\left( \Gamma \bar{\gamma}\right) \func{Re}\left( \Phi \bar{%
\phi}\right) }}\left\vert \left\langle x,e\right\rangle \left\langle
e,y\right\rangle \right\vert .  \label{1.13}
\end{equation}%
Here the constant $\frac{1}{4}$ is also best possible.

In the case that both $\left\langle x,e\right\rangle ,$ $\left\langle
e,y\right\rangle \neq 0$ (which is actually the interesting case), we have%
\begin{equation}
\left\vert \frac{\left\langle x,y\right\rangle }{\left\langle
x,e\right\rangle \left\langle e,y\right\rangle }-1\right\vert \leq \frac{1}{4%
}\cdot \frac{\left\vert \Gamma -\gamma \right\vert \left\vert \Phi -\phi
\right\vert }{\sqrt{\func{Re}\left( \Gamma \bar{\gamma}\right) \func{Re}%
\left( \Phi \bar{\phi}\right) }}.  \label{1.14}
\end{equation}

Now, for an orthornormal family of vectors in $H,$ i.e., we recall that $%
\left\langle e_{i},e_{j}\right\rangle =0$ if $i,j\in \mathbb{N}$, $i\neq j$
and $\left\Vert e_{i}\right\Vert =1$ for $i\in \mathbb{N}$, the following
inequality, called the \textit{Bessel inequality}%
\begin{equation}
\sum_{i=1}^{\infty }\left\vert \left\langle x,e_{i}\right\rangle \right\vert
^{2}\leq \left\Vert x\right\Vert ^{2},\ \ \ \ \ \ x\in H;  \label{1.15}
\end{equation}%
holds.

If $\left( H;\left\langle \cdot ,\cdot \right\rangle \right) $ is an
infinite dimensional Hilbert space over the real or complex number field $%
\mathbb{K}$, $\left( e_{i}\right) _{i\in \mathbb{N}}$ an orthornormal family
in $H,$ $\mathbf{\lambda }=\left( \lambda _{i}\right) _{i\in \mathbb{N}}\in
\ell ^{2}\left( \mathbb{K}\right) $ and $r>0$ with the property that%
\begin{equation}
\sum_{i=1}^{\infty }\left\vert \lambda _{i}\right\vert ^{2}>r^{2},
\label{1.16}
\end{equation}%
then, for $x\in H$ such that%
\begin{equation}
\left\Vert x-\sum_{i=1}^{\infty }\lambda _{i}e_{i}\right\Vert \leq r,
\label{1.17}
\end{equation}%
holds, one has the inequalities \cite{SSDa}%
\begin{align}
\left\Vert x\right\Vert ^{2}& \leq \frac{\left( \sum_{i=1}^{\infty }\func{Re}%
\left[ \bar{\lambda}_{i}\left\langle x,e_{i}\right\rangle \right] \right)
^{2}}{\sum_{i=1}^{\infty }\left\vert \lambda _{i}\right\vert ^{2}-r^{2}}\leq 
\frac{\left\vert \sum_{i=1}^{\infty }\bar{\lambda}_{i}\left\langle
x,e_{i}\right\rangle \right\vert ^{2}}{\sum_{i=1}^{\infty }\left\vert
\lambda _{i}\right\vert ^{2}-r^{2}}  \label{1.18} \\
&  \notag \\
& \leq \frac{\sum_{i=1}^{\infty }\left\vert \lambda _{i}\right\vert ^{2}}{%
\sum_{i=1}^{\infty }\left\vert \lambda _{i}\right\vert ^{2}-r^{2}}%
\sum_{i=1}^{\infty }\left\vert \left\langle x,e_{i}\right\rangle \right\vert
^{2}.  \notag
\end{align}%
An additive version of interest is \cite{SSDa}%
\begin{equation}
0\leq \left\Vert x\right\Vert ^{2}-\sum_{i=1}^{\infty }\left\vert
\left\langle x,e_{i}\right\rangle \right\vert ^{2}\leq \frac{r^{2}}{%
\sum_{i=1}^{\infty }\left\vert \lambda _{i}\right\vert ^{2}-r^{2}}%
\sum_{i=1}^{\infty }\left\vert \left\langle x,e_{i}\right\rangle \right\vert
^{2}.  \label{1.17b}
\end{equation}%
Finally, if $\mathbf{\Gamma }=\left( \Gamma _{i}\right) _{i\in \mathbb{N}},$ 
$\mathbf{\gamma }=\left( \gamma _{i}\right) _{i\in \mathbb{N}}\in \ell
^{2}\left( \mathbb{K}\right) $ are such that $\sum_{i=1}^{\infty }\func{Re}%
\left( \Gamma _{i}\overline{\gamma _{i}}\right) >0$ and for $x\in H,$ either%
\begin{equation}
\left\Vert x-\sum_{i=1}^{\infty }\frac{\gamma _{i}+\Gamma _{i}}{2}%
e_{i}\right\Vert \leq \frac{1}{2}\left( \sum_{i=1}^{\infty }\left\vert
\Gamma _{i}-\gamma _{i}\right\vert ^{2}\right) ^{\frac{1}{2}}  \label{1.18b}
\end{equation}%
or, equivalently,%
\begin{equation}
\func{Re}\left\langle \sum_{i=1}^{\infty }\Gamma
_{i}e_{i}-x,x-\sum_{i=1}^{\infty }\gamma _{i}e_{i}\right\rangle \geq 0
\label{1.19}
\end{equation}%
holds, then \cite{SSDa}%
\begin{align}
\left\Vert x\right\Vert ^{2}& \leq \frac{1}{4}\cdot \frac{\left(
\sum_{i=1}^{\infty }\func{Re}\left[ \left( \bar{\Gamma}_{i}+\bar{\gamma}%
_{i}\right) \left\langle x,e_{i}\right\rangle \right] \right) ^{2}}{%
\sum_{i=1}^{\infty }\func{Re}\left( \Gamma _{i}\bar{\gamma}_{i}\right) }\leq 
\frac{1}{4}\cdot \frac{\left\vert \sum_{i=1}^{\infty }\left( \bar{\Gamma}%
_{i}+\bar{\gamma}_{i}\right) \left\langle x,e_{i}\right\rangle \right\vert
^{2}}{\sum_{i=1}^{\infty }\func{Re}\left( \Gamma _{i}\bar{\gamma}_{i}\right) 
}  \label{1.20} \\
&  \notag \\
& \leq \frac{1}{4}\cdot \frac{\sum_{i=1}^{\infty }\left\vert \Gamma
_{i}+\gamma _{i}\right\vert ^{2}}{\sum_{i=1}^{\infty }\func{Re}\left( \Gamma
_{i}\bar{\gamma}_{i}\right) }\sum_{i=1}^{\infty }\left\vert \left\langle
x,e_{i}\right\rangle \right\vert ^{2}.  \notag
\end{align}%
The constant $\frac{1}{4}$ is best possible in all inequalities (\ref{1.20}).

The following additive version may be stated as well \cite{SSDa}%
\begin{equation}
0\leq \left\Vert x\right\Vert ^{2}-\sum_{i=1}^{\infty }\left\vert
\left\langle x,e_{i}\right\rangle \right\vert ^{2}\leq \frac{1}{4}\cdot 
\frac{\sum_{i=1}^{\infty }\left\vert \Gamma _{i}-\gamma _{i}\right\vert ^{2}%
}{\sum_{i=1}^{\infty }\func{Re}\left( \Gamma _{i}\bar{\gamma}_{i}\right) }%
\sum_{i=1}^{\infty }\left\vert \left\langle x,e_{i}\right\rangle \right\vert
^{2}.  \label{1.21}
\end{equation}%
Here the constant $\frac{1}{4}$ is also best possible.

The present paper is a continuation of \cite{SSDa}. Here we point out
different reverses of the Schwarz, triangle and Bessel inequalities that are
also sharp. Applications for Gr\"{u}ss type inequalities are provided. Some
integral inequalities that are natural consequences of the above, are stated
as well.

\section{New Reverses of Schwarz's Inequality}

The following result holds.

\begin{theorem}
\label{t2.1}Let $\left( H;\left\langle \cdot ,\cdot \right\rangle \right) $
be an inner product space over the real or complex number field $\mathbb{K}$%
, $x,a\in H$ and $r>0.$ If%
\begin{equation}
x\in \bar{B}\left( a,r\right) :=\left\{ z\in H|\left\Vert z-a\right\Vert
\leq r\right\} ,  \label{2.1}
\end{equation}%
then we have the inequalities:%
\begin{align}
0& \leq \left\Vert x\right\Vert \left\Vert a\right\Vert -\left\vert
\left\langle x,a\right\rangle \right\vert \leq \left\Vert x\right\Vert
\left\Vert a\right\Vert -\left\vert \func{Re}\left\langle x,a\right\rangle
\right\vert  \label{2.2} \\
& \leq \left\Vert x\right\Vert \left\Vert a\right\Vert -\func{Re}%
\left\langle x,a\right\rangle \leq \frac{1}{2}r^{2}.  \notag
\end{align}%
The constant $\frac{1}{2}$ is best possible in (\ref{2.2}) in the sense that
it cannot be replaced by a smaller constant.
\end{theorem}

\begin{proof}
The condition (\ref{2.1}) is clearly equivalent to%
\begin{equation}
\left\Vert x\right\Vert ^{2}+\left\Vert a\right\Vert ^{2}\leq 2\func{Re}%
\left\langle x,a\right\rangle +r^{2}.  \label{2.3}
\end{equation}%
Using the elementary inequality%
\begin{equation}
2\left\Vert x\right\Vert \left\Vert a\right\Vert \leq \left\Vert
x\right\Vert ^{2}+\left\Vert a\right\Vert ^{2},\ \ \ \ a,x\in H  \label{2.4}
\end{equation}%
and (\ref{2.3}), we deduce%
\begin{equation}
2\left\Vert x\right\Vert \left\Vert a\right\Vert \leq 2\func{Re}\left\langle
x,a\right\rangle +r^{2},  \label{2.4a}
\end{equation}%
giving the last inequality in (\ref{2.2}). The other inequalities are
obvious.

To prove the sharpness of the constant $\frac{1}{2},$ assume that%
\begin{equation}
0\leq \left\Vert x\right\Vert \left\Vert a\right\Vert -\func{Re}\left\langle
x,a\right\rangle \leq cr^{2}  \label{2.3a}
\end{equation}%
for any $x,a\in H$ and $r>0$ satisfying (\ref{2.1}).

Assume that $a,e\in H,$ $\left\Vert a\right\Vert =\left\Vert e\right\Vert =1$
and $e\perp a.$ If $r=\sqrt{\varepsilon },$ $\varepsilon >0$ and if we
define $x=a+\sqrt{\varepsilon }e,$ then $\left\Vert x-a\right\Vert =\sqrt{%
\varepsilon }=r$ showing that the condition (\ref{2.1}) is fulfilled.

On the other hand,%
\begin{align*}
\left\Vert x\right\Vert \left\Vert a\right\Vert -\func{Re}\left\langle
x,a\right\rangle & =\sqrt{\left\Vert a+\sqrt{\varepsilon }e\right\Vert ^{2}}-%
\func{Re}\left\langle a+\sqrt{\varepsilon }e,a\right\rangle \\
& =\sqrt{\left\Vert a\right\Vert ^{2}+\varepsilon \left\Vert e\right\Vert
^{2}}-\left\Vert a\right\Vert ^{2} \\
& =\sqrt{1+\varepsilon }-1.
\end{align*}%
Utilising (\ref{2.3a}), we conclude that%
\begin{equation}
\sqrt{1+\varepsilon }-1\leq c\varepsilon \text{ \ for any }\varepsilon >0.
\label{2.4b}
\end{equation}%
Multiplying (\ref{2.4b}) by $\sqrt{1+\varepsilon }+1>0$ and then dividing by 
$\varepsilon >0,$ we get%
\begin{equation}
\left( \sqrt{1+\varepsilon }+1\right) c\geq 1\text{ \ for any \ }\varepsilon
>0.  \label{2.5}
\end{equation}%
Letting $\varepsilon \rightarrow 0+$ in (\ref{2.5}), we deduce $c\geq \frac{1%
}{2},$ and the theorem is proved.
\end{proof}

The following result also holds.

\begin{theorem}
\label{t2.2}Let $\left( H;\left\langle \cdot ,\cdot \right\rangle \right) $
be an inner product space over $\mathbb{K}$ and $x,y\in H,$ $\gamma ,\Gamma
\in \mathbb{K}$ \ $\left( \Gamma \neq \gamma ,-\gamma \right) $ so that
either%
\begin{equation}
\func{Re}\left\langle \Gamma y-x,x-\gamma y\right\rangle \geq 0,  \label{2.6}
\end{equation}%
or, equivalently,%
\begin{equation}
\left\Vert x-\frac{\gamma +\Gamma }{2}y\right\Vert \leq \frac{1}{2}%
\left\vert \Gamma -\gamma \right\vert \left\Vert y\right\Vert  \label{2.7}
\end{equation}%
holds. Then we have the inequalities%
\begin{align}
0& \leq \left\Vert x\right\Vert \left\Vert y\right\Vert -\left\vert
\left\langle x,y\right\rangle \right\vert  \label{2.8} \\
& \leq \left\Vert x\right\Vert \left\Vert y\right\Vert -\left\vert \func{Re}%
\left[ \frac{\bar{\Gamma}+\bar{\gamma}}{\left\vert \Gamma +\gamma
\right\vert }\left\langle x,y\right\rangle \right] \right\vert  \notag \\
& \leq \left\Vert x\right\Vert \left\Vert y\right\Vert -\func{Re}\left[ 
\frac{\bar{\Gamma}+\bar{\gamma}}{\left\vert \Gamma +\gamma \right\vert }%
\left\langle x,y\right\rangle \right]  \notag \\
& \leq \frac{1}{4}\cdot \frac{\left\vert \Gamma -\gamma \right\vert ^{2}}{%
\left\vert \Gamma +\gamma \right\vert }\left\Vert y\right\Vert ^{2}.  \notag
\end{align}%
The constant $\frac{1}{4}$ in the last inequality is best possible.
\end{theorem}

\begin{proof}
The proof of the equivalence between the inequalities (\ref{2.6}) and (\ref%
{2.7}) follows by the fact that in an inner product space $\func{Re}%
\left\langle Z-x,x-z\right\rangle \geq 0$ for $x,z,Z\in H$ is equivalent to 
\begin{equation*}
\left\Vert x-\frac{z+Z}{2}\right\Vert \leq \frac{1}{2}\left\Vert
Z-z\right\Vert ,
\end{equation*}%
(see for example \cite{SSD3}).

Consider for $a,y\neq 0,$ $a=\frac{\Gamma +\gamma }{2}\cdot y$ and $r=\frac{1%
}{2}\left\vert \Gamma -\gamma \right\vert \left\Vert y\right\Vert .$ Thus
from (\ref{2.2}), we get%
\begin{align*}
0& \leq \left\Vert x\right\Vert \left\vert \frac{\Gamma +\gamma }{2}%
\right\vert \left\Vert y\right\Vert -\left\vert \frac{\Gamma +\gamma }{2}%
\right\vert \left\vert \left\langle x,y\right\rangle \right\vert \\
& \leq \left\Vert x\right\Vert \left\vert \frac{\Gamma +\gamma }{2}%
\right\vert \left\Vert y\right\Vert -\left\vert \func{Re}\left[ \frac{\bar{%
\Gamma}+\bar{\gamma}}{\left\vert \Gamma +\gamma \right\vert }\left\langle
x,y\right\rangle \right] \right\vert \\
& \leq \left\Vert x\right\Vert \left\vert \frac{\Gamma +\gamma }{2}%
\right\vert \left\Vert y\right\Vert -\func{Re}\left[ \frac{\bar{\Gamma}+\bar{%
\gamma}}{\left\vert \Gamma +\gamma \right\vert }\left\langle
x,y\right\rangle \right] \\
& \leq \frac{1}{8}\cdot \left\vert \Gamma -\gamma \right\vert ^{2}\left\Vert
y\right\Vert ^{2}.
\end{align*}%
Dividing by $\frac{1}{2}\left\vert \Gamma +\gamma \right\vert \geq 0,$ we
deduce the desired inequality (\ref{2.8}).

To prove the sharpness of the constant $\frac{1}{4},$ assume that there
exists a $c>0$ such that:%
\begin{equation}
\left\Vert x\right\Vert \left\Vert y\right\Vert -\func{Re}\left[ \frac{\bar{%
\Gamma}+\bar{\gamma}}{\left\vert \Gamma +\gamma \right\vert }\left\langle
x,y\right\rangle \right] \leq c\cdot \frac{\left\vert \Gamma -\gamma
\right\vert ^{2}}{\left\vert \Gamma +\gamma \right\vert }\left\Vert
y\right\Vert ^{2},  \label{2.9}
\end{equation}%
provided either (\ref{2.6}) or (\ref{2.7}) holds.

Consider the real inner product space $\left( \mathbb{R}^{2},\left\langle
\cdot ,\cdot \right\rangle \right) $ with $\left\langle \mathbf{\bar{x}},%
\mathbf{\bar{y}}\right\rangle =x_{1}y_{1}+x_{2}y_{2},$ $\mathbf{\bar{x}}%
=\left( x_{1},x_{2}\right) ,$ $\mathbf{\bar{y}}=\left( y_{1},y_{2}\right)
\in \mathbb{R}^{2}.$ Let $y=\left( 1,1\right) $ and $\Gamma ,\gamma >0$ with 
$\Gamma >\gamma .$ Then, by (\ref{2.9}), we deduce%
\begin{equation}
\sqrt{2}\sqrt{x_{1}^{2}+x_{2}^{2}}-\left( x_{1}+x_{2}\right) \leq 2c\cdot 
\frac{\left( \Gamma -\gamma \right) ^{2}}{\Gamma +\gamma }.  \label{2.10}
\end{equation}%
If $x_{1}=\Gamma ,$ $x_{2}=\gamma ,$ then 
\begin{equation*}
\left\langle \Gamma \mathbf{\bar{y}}-\mathbf{\bar{x}},\mathbf{\bar{x}}%
-\gamma \mathbf{\bar{y}}\right\rangle =\left( \Gamma -x_{1}\right) \left(
x_{1}-\gamma \right) +\left( \Gamma -x_{2}\right) \left( x_{2}-\gamma
\right) =0,
\end{equation*}%
showing that the condition (\ref{2.6}) is valid. Replacing \thinspace $x_{1}$
and $x_{2}$ in (\ref{2.10}), we deduce%
\begin{equation}
\sqrt{2}\sqrt{\Gamma ^{2}+\gamma ^{2}}-\left( \Gamma +\gamma \right) \leq 2c%
\frac{\left( \Gamma -\gamma \right) ^{2}}{\Gamma +\gamma }.  \label{2.11}
\end{equation}%
If in (\ref{2.11}) we choose $\Gamma =1+\varepsilon ,$ $\gamma
=1-\varepsilon $ with $\varepsilon \in \left( 0,1\right) ,$ we deduce 
\begin{equation*}
2\sqrt{1+\varepsilon ^{2}}-2\leq 2c\frac{4\varepsilon ^{2}}{2},
\end{equation*}%
giving%
\begin{equation}
\sqrt{1+\varepsilon ^{2}}-1\leq 2c\varepsilon ^{2}.  \label{2.12}
\end{equation}%
Finally, multiplying (\ref{2.12}) with $\sqrt{1+\varepsilon ^{2}}+1>0$ and
thus dividing by $\varepsilon ^{2},$ we deduce%
\begin{equation}
1\leq 2c\left( \sqrt{1+\varepsilon ^{2}}+1\right) \text{ \ for any \ }%
\varepsilon \in \left( 0,1\right) .  \label{2.13}
\end{equation}%
Letting $\varepsilon \rightarrow 0+$ in (\ref{2.13}) we get $c\geq \frac{1}{4%
},$ and the sharpness of the constant is proved.
\end{proof}

For some recent results in connection to Schwarz's inequality, see \cite{ADR}%
, \cite{DM} and \cite{GH}.

\section{Reverses of the Triangle Inequality}

The following reverse of the triangle inequality in inner product spaces
holds.

\begin{proposition}
\label{p2.3}Let $\left( H;\left\langle \cdot ,\cdot \right\rangle \right) $
be an inner product space over the real or complex number field $\mathbb{K}$%
, $x,a\in H$ and $r>0.$ If $\left\Vert x-a\right\Vert \leq r,$ then we have
the inequality%
\begin{equation}
0\leq \left\Vert x\right\Vert +\left\Vert a\right\Vert -\left\Vert
x+a\right\Vert \leq r.  \label{2.14}
\end{equation}
\end{proposition}

\begin{proof}
Since%
\begin{equation*}
\left( \left\Vert x\right\Vert +\left\Vert a\right\Vert \right)
^{2}-\left\Vert x+a\right\Vert ^{2}\leq 2\left( \left\Vert x\right\Vert
\left\Vert a\right\Vert -\func{Re}\left\langle x,a\right\rangle \right) ,
\end{equation*}%
then by Theorem \ref{t2.1} we deduce%
\begin{equation}
\left( \left\Vert x\right\Vert +\left\Vert a\right\Vert \right)
^{2}-\left\Vert x+a\right\Vert ^{2}\leq r^{2},  \label{2.15}
\end{equation}%
from where we obtain%
\begin{equation}
\left\Vert x\right\Vert +\left\Vert a\right\Vert \leq \sqrt{r^{2}+\left\Vert
x+a\right\Vert ^{2}}\leq r+\left\Vert x+a\right\Vert ,  \label{2.16}
\end{equation}%
giving the desired result (\ref{2.14}).
\end{proof}

We may state the following result.

\begin{proposition}
\label{p2.4}Let $\left( H;\left\langle \cdot ,\cdot \right\rangle \right) $
be an inner product space over $\mathbb{K}$ and $x,y\in H,$ $M>m>0$ such
that either%
\begin{equation}
\func{Re}\left\langle My-x,x-my\right\rangle \geq 0,  \label{2.17}
\end{equation}%
or, equivalently,%
\begin{equation}
\left\Vert x-\frac{M+m}{2}y\right\Vert \leq \frac{1}{2}\left( M-m\right)
\left\Vert y\right\Vert ,  \label{2.18}
\end{equation}%
holds. Then we have the inequality%
\begin{equation}
0\leq \left\Vert x\right\Vert +\left\Vert y\right\Vert -\left\Vert
x+y\right\Vert \leq \frac{\sqrt{2}}{2}\cdot \frac{\left( M-m\right) }{\sqrt{%
M+m}}\left\Vert y\right\Vert .  \label{2.19}
\end{equation}
\end{proposition}

\begin{proof}
By Theorem \ref{t2.2} for $\Gamma =M,$ $\gamma =m,$ we have the inequality%
\begin{equation}
\left\Vert x\right\Vert \left\Vert y\right\Vert -\func{Re}\left\langle
x,y\right\rangle \leq \frac{1}{4}\cdot \frac{\left( M-m\right) ^{2}}{\left(
M+m\right) }\left\Vert y\right\Vert ^{2}.  \label{2.20}
\end{equation}%
Then we may state that%
\begin{align*}
\left( \left\Vert x\right\Vert +\left\Vert y\right\Vert \right)
^{2}-\left\Vert x+y\right\Vert ^{2}& =2\left( \left\Vert x\right\Vert
\left\Vert y\right\Vert -\func{Re}\left\langle x,y\right\rangle \right) \\
& \leq \frac{1}{2}\cdot \frac{\left( M-m\right) ^{2}}{M+m}\left\Vert
y\right\Vert ^{2}
\end{align*}%
from where we get%
\begin{align}
\left\Vert x\right\Vert +\left\Vert y\right\Vert & \leq \sqrt{\frac{1}{2}%
\cdot \frac{\left( M-m\right) ^{2}}{M+m}\left\Vert y\right\Vert
^{2}+\left\Vert x+y\right\Vert ^{2}}  \label{2.21} \\
& \leq \left\Vert x+y\right\Vert +\frac{\left( M-m\right) }{\sqrt{2\left(
M+m\right) }}\left\Vert y\right\Vert ,  \notag
\end{align}%
giving the desired inequality (\ref{2.19}).
\end{proof}

For some results related to triangle inequality in inner product spaces, see 
\cite{JBDFTM}, \cite{SMK}, \cite{PMM} and \cite{DKR}.

\section{Some Gr\"{u}ss Type Inequalities}

We may state the following result.

\begin{theorem}
\label{t4.1}Let $\left( H;\left\langle \cdot ,\cdot \right\rangle \right) $
be an inner product space over the real or complex number field $\mathbb{K}$
and $x,y,e\in H$ with $\left\Vert e\right\Vert =1.$ If $r_{1},r_{2}>0$ and 
\begin{equation}
\left\Vert x-e\right\Vert \leq r_{1},\ \ \ \ \left\Vert y-e\right\Vert \leq
r_{2},  \label{4.1}
\end{equation}%
then we have the inequalities%
\begin{eqnarray}
\left\vert \left\langle x,y\right\rangle -\left\langle x,e\right\rangle
\left\langle e,y\right\rangle \right\vert &\leq &\frac{1}{2}r_{1}r_{2}\sqrt{%
\left\Vert x\right\Vert +\left\vert \left\langle x,e\right\rangle
\right\vert }\cdot \sqrt{\left\Vert y\right\Vert +\left\vert \left\langle
y,e\right\rangle \right\vert }  \label{4.2} \\
&\leq &r_{1}r_{2}\left\Vert x\right\Vert \left\Vert y\right\Vert .  \notag
\end{eqnarray}%
The constant $\frac{1}{2}$ is best possible in the sense that it cannot be
replaced by a smaller constant.
\end{theorem}

\begin{proof}
Apply Schwarz's inequality for the vectors $x-\left\langle x,e\right\rangle
e,$ $y-\left\langle y,e\right\rangle e$ to get (see also \cite{SSD3}):%
\begin{equation}
\left\vert \left\langle x,y\right\rangle -\left\langle x,e\right\rangle
\left\langle e,y\right\rangle \right\vert ^{2}\leq \left( \left\Vert
x\right\Vert ^{2}-\left\vert \left\langle x,e\right\rangle \right\vert
^{2}\right) \left( \left\Vert y\right\Vert ^{2}-\left\vert \left\langle
y,e\right\rangle \right\vert ^{2}\right) .  \label{4.3}
\end{equation}%
Using Theorem \ref{t2.1} for $a=e,$ we have%
\begin{align}
0& \leq \left\Vert x\right\Vert ^{2}-\left\vert \left\langle
x,e\right\rangle \right\vert ^{2}  \label{4.4} \\
& =\left( \left\Vert x\right\Vert -\left\vert \left\langle x,e\right\rangle
\right\vert \right) \left( \left\Vert x\right\Vert +\left\vert \left\langle
x,e\right\rangle \right\vert \right)  \notag \\
& \leq \frac{1}{2}r_{1}^{2}\left( \left\Vert x\right\Vert +\left\vert
\left\langle x,e\right\rangle \right\vert \right) \leq r_{1}^{2}\left\Vert
x\right\Vert ,  \notag
\end{align}%
and, in a similar way%
\begin{align}
0& \leq \left\Vert y\right\Vert ^{2}-\left\vert \left\langle
y,e\right\rangle \right\vert ^{2}  \label{4.5} \\
& \leq \frac{1}{2}r_{2}^{2}\left( \left\Vert y\right\Vert +\left\vert
\left\langle y,e\right\rangle \right\vert \right) \leq r_{2}^{2}\left\Vert
y\right\Vert .  \notag
\end{align}%
Utilising (\ref{4.3}) -- (\ref{4.5}), we may state that%
\begin{align}
\left\vert \left\langle x,y\right\rangle -\left\langle x,e\right\rangle
\left\langle e,y\right\rangle \right\vert ^{2}& \leq \frac{1}{4}%
r_{1}^{2}r_{2}^{2}\left( \left\Vert x\right\Vert +\left\vert \left\langle
x,e\right\rangle \right\vert \right) \left( \left\Vert y\right\Vert
+\left\vert \left\langle y,e\right\rangle \right\vert \right)  \label{4.6} \\
& \leq r_{1}^{2}r_{2}^{2}\left\Vert x\right\Vert \left\Vert y\right\Vert , 
\notag
\end{align}%
giving the desired inequality (\ref{4.2}).

To prove the sharpness of the constant $\frac{1}{2}$, let assume $x=y$ in (%
\ref{4.2}), to get%
\begin{equation}
\left\Vert x\right\Vert ^{2}-\left\vert \left\langle x,e\right\rangle
\right\vert ^{2}\leq \frac{1}{2}r_{1}^{2}\left( \left\Vert x\right\Vert
+\left\vert \left\langle x,e\right\rangle \right\vert \right) ,  \label{4.7}
\end{equation}%
provided $\left\Vert x-e\right\Vert \leq r_{1}.$ If $x\neq 0,$ then dividing
(\ref{4.7}) with $\left\Vert x\right\Vert +\left\vert \left\langle
x,e\right\rangle \right\vert >0$ we get%
\begin{equation}
\left\Vert x\right\Vert -\left\vert \left\langle x,e\right\rangle
\right\vert \leq \frac{1}{2}r_{1}^{2}  \label{4.8}
\end{equation}%
provided $\left\Vert x-e\right\Vert \leq r_{1},$ $\left\Vert e\right\Vert
=1. $ However, (\ref{4.8}) is in fact (\ref{2.2}) for $a=e,$ for which we
have shown that $\frac{1}{2}$ is the best possible constant.
\end{proof}

The following result also holds.

\begin{theorem}
\label{t4.2}With the assumptions of Theorem \ref{t4.1}, we have the
inequality%
\begin{equation}
\left\vert \left\langle x,y\right\rangle -\left\langle x,e\right\rangle
\left\langle e,y\right\rangle \right\vert \leq r_{1}r_{2}\sqrt{\frac{1}{4}%
r_{1}^{2}+\left\vert \left\langle x,e\right\rangle \right\vert }\cdot \sqrt{%
\frac{1}{4}r_{2}^{2}+\left\vert \left\langle y,e\right\rangle \right\vert }.
\label{4.9}
\end{equation}
\end{theorem}

\begin{proof}
Note that, from Theorem \ref{t2.2}, we have%
\begin{equation}
\left\Vert x\right\Vert \left\Vert a\right\Vert \leq \left\vert \left\langle
x,a\right\rangle \right\vert +\frac{1}{2}r^{2}  \label{4.10}
\end{equation}%
provided $\left\Vert x-a\right\Vert \leq r.$

Taking the square in (\ref{4.10}) and arranging the terms, we obtain:%
\begin{equation}
0\leq \left\Vert x\right\Vert ^{2}\left\Vert a\right\Vert ^{2}-\left\vert
\left\langle x,a\right\rangle \right\vert ^{2}\leq r^{2}\left( \frac{1}{4}%
r^{2}+\left\vert \left\langle x,a\right\rangle \right\vert \right) ,
\label{4.11}
\end{equation}%
provided $\left\Vert x-a\right\Vert \leq r.$

Using the assumption of the theorem, we then have%
\begin{equation}
0\leq \left\Vert x\right\Vert ^{2}-\left\vert \left\langle x,e\right\rangle
\right\vert ^{2}\leq r_{1}^{2}\left( \frac{1}{4}r_{1}^{2}+\left\vert
\left\langle x,e\right\rangle \right\vert \right) ,  \label{4.12}
\end{equation}%
and 
\begin{equation}
0\leq \left\Vert y\right\Vert ^{2}-\left\vert \left\langle y,e\right\rangle
\right\vert ^{2}\leq r_{2}^{2}\left( \frac{1}{4}r_{2}^{2}+\left\vert
\left\langle y,e\right\rangle \right\vert \right) .  \label{4.13}
\end{equation}%
Utilising (\ref{4.3}), (\ref{4.12}) and (\ref{4.13}), we deduce the desired
inequality (\ref{4.9}).
\end{proof}

The following result may be stated as well.

\begin{theorem}
\label{t4.3}Let $\left( H;\left\langle \cdot ,\cdot \right\rangle \right) $
be an inner product space over $\mathbb{K}$ and $x,y,e\in H$ with $%
\left\Vert e\right\Vert =1.$ Suppose also that $a,A,b,B\in \mathbb{K}$ $%
\left( \mathbb{K}=\mathbb{C},\mathbb{R}\right) $ so that $A\neq \pm a,B\neq
\pm b.$ If either%
\begin{equation}
\func{Re}\left\langle Ae-x,x-ae\right\rangle \geq 0,\ \ \func{Re}%
\left\langle Be-y,y-be\right\rangle \geq 0,  \label{4.14}
\end{equation}%
or, equivalently,%
\begin{equation}
\left\Vert x-\frac{a+A}{2}e\right\Vert \leq \frac{1}{2}\left\vert
A-a\right\vert ,\ \ \ \left\Vert y-\frac{b+B}{2}e\right\Vert \leq \frac{1}{2}%
\left\vert B-b\right\vert ,  \label{4.15}
\end{equation}%
holds, then we have the inequality%
\begin{align}
& \left\vert \left\langle x,y\right\rangle -\left\langle x,e\right\rangle
\left\langle e,y\right\rangle \right\vert  \label{4.16} \\
& \leq \frac{1}{4}\cdot \frac{\left\vert A-a\right\vert \left\vert
B-b\right\vert }{\sqrt{\left\vert A+a\right\vert \left\vert B+b\right\vert }}%
\sqrt{\left\Vert x\right\Vert +\left\vert \left\langle x,e\right\rangle
\right\vert }\cdot \sqrt{\left\Vert y\right\Vert +\left\vert \left\langle
y,e\right\rangle \right\vert }  \notag \\
& \leq \frac{1}{2}\cdot \frac{\left\vert A-a\right\vert \left\vert
B-b\right\vert }{\sqrt{\left\vert A+a\right\vert \left\vert B+b\right\vert }}%
\sqrt{\left\Vert x\right\Vert \left\Vert y\right\Vert }.  \notag
\end{align}%
The constant $\frac{1}{4}$ is best possible in (\ref{4.16}).
\end{theorem}

\begin{proof}
From Theorem \ref{t2.2}, we may state that%
\begin{align}
0& \leq \left\Vert x\right\Vert ^{2}-\left\vert \left\langle
x,e\right\rangle \right\vert ^{2}  \label{4.17} \\
& =\left( \left\Vert x\right\Vert -\left\vert \left\langle x,e\right\rangle
\right\vert \right) \left( \left\Vert x\right\Vert +\left\vert \left\langle
x,e\right\rangle \right\vert \right)  \notag \\
& \leq \frac{1}{4}\cdot \frac{\left\vert A-a\right\vert ^{2}}{\left\vert
A+a\right\vert }\left( \left\Vert x\right\Vert +\left\vert \left\langle
x,e\right\rangle \right\vert \right) ,  \notag
\end{align}%
and%
\begin{equation}
0\leq \left\Vert y\right\Vert ^{2}-\left\vert \left\langle y,e\right\rangle
\right\vert ^{2}\leq \frac{1}{4}\cdot \frac{\left\vert B-b\right\vert ^{2}}{%
\left\vert B+b\right\vert }\left( \left\Vert y\right\Vert +\left\vert
\left\langle y,e\right\rangle \right\vert \right) .  \label{4.18}
\end{equation}%
Making use of (\ref{4.3}) and (\ref{4.17}), (\ref{4.18}), we deduce the
first inequality in (\ref{4.16}).

The best constant follows by the use of Theorem \ref{t2.2}, and we omit the
details.
\end{proof}

Finally, we may state the following theorem as well.

\begin{theorem}
\label{t4.4}With the assumptions of Theorem \ref{t4.3}, we have the
inequality%
\begin{multline}
\left\vert \left\langle x,y\right\rangle -\left\langle x,e\right\rangle
\left\langle e,y\right\rangle \right\vert  \label{4.19} \\
\leq \frac{1}{2}\cdot \frac{\left\vert A-a\right\vert \left\vert
B-b\right\vert }{\sqrt{\left\vert A+a\right\vert \left\vert B+b\right\vert }}%
\sqrt{\frac{1}{8}\cdot \frac{\left\vert A-a\right\vert ^{2}}{\left\vert
A+a\right\vert }+\left\vert \left\langle x,e\right\rangle \right\vert }\cdot 
\sqrt{\frac{1}{8}\cdot \frac{\left\vert B-b\right\vert ^{2}}{\left\vert
B+b\right\vert }+\left\vert \left\langle y,e\right\rangle \right\vert }.
\end{multline}
\end{theorem}

\begin{proof}
Using Theorem \ref{t2.1}, we may state that%
\begin{equation*}
0\leq \left\Vert x\right\Vert -\left\vert \left\langle x,e\right\rangle
\right\vert \leq \frac{1}{4}\cdot \frac{\left\vert A-a\right\vert ^{2}}{%
\left\vert A+a\right\vert }.
\end{equation*}%
This inequality implies that%
\begin{equation*}
\left\Vert x\right\Vert ^{2}\leq \left\vert \left\langle x,e\right\rangle
\right\vert ^{2}+\frac{1}{2}\left\vert \left\langle x,e\right\rangle
\right\vert \cdot \frac{\left\vert A-a\right\vert ^{2}}{\left\vert
A+a\right\vert }+\frac{1}{16}\cdot \frac{\left\vert A-a\right\vert ^{4}}{%
\left\vert A+a\right\vert ^{2}}
\end{equation*}%
giving%
\begin{equation}
0\leq \left\Vert x\right\Vert ^{2}-\left\vert \left\langle x,e\right\rangle
\right\vert ^{2}\leq \frac{1}{2}\cdot \frac{\left\vert A-a\right\vert ^{2}}{%
\left\vert A+a\right\vert }\left[ \left\vert \left\langle x,e\right\rangle
\right\vert +\frac{1}{8}\cdot \frac{\left\vert A-a\right\vert ^{2}}{%
\left\vert A+a\right\vert }\right] .  \label{4.20}
\end{equation}%
Similarly, we have%
\begin{equation}
0\leq \left\Vert y\right\Vert ^{2}-\left\vert \left\langle y,e\right\rangle
\right\vert ^{2}\leq \frac{1}{2}\cdot \frac{\left\vert B-b\right\vert ^{2}}{%
\left\vert B+b\right\vert }\left[ \left\vert \left\langle y,e\right\rangle
\right\vert +\frac{1}{8}\cdot \frac{\left\vert B-b\right\vert ^{2}}{%
\left\vert B+b\right\vert }\right] .  \label{4.21}
\end{equation}%
By making use of (\ref{4.3}) and (\ref{4.20}), (\ref{4.21}), we deduce the
desired inequality (\ref{4.19}).
\end{proof}

For some recent results on Gr\"{u}ss type inequalities in inner product
spaces, see \cite{SSD0}, \cite{SSD00} and \cite{PFR}.

\section{Reverses of Bessel's Inequality}

Let $\left( H;\left\langle \cdot ,\cdot \right\rangle \right) $ be a real or
complex infinite dimensional Hilbert space and $\left( e_{i}\right) _{i\in 
\mathbb{N}}$ an orthornormal family in $H,$ i.e., we recall that $%
\left\langle e_{i},e_{j}\right\rangle =0$ if $i,j\in \mathbb{N}$, $i\neq j$
and $\left\Vert e_{i}\right\Vert =1$ for $i\in \mathbb{N}$.

It is well known that, if $x\in H,$ then the series $\sum_{i=1}^{\infty
}\left\vert \left\langle x,e_{i}\right\rangle \right\vert ^{2}$ is
convergent and the following inequality, called \textit{Bessel's inequality, 
}%
\begin{equation}
\sum_{i=1}^{\infty }\left\vert \left\langle x,e_{i}\right\rangle \right\vert
^{2}\leq \left\Vert x\right\Vert ^{2},  \label{5.1}
\end{equation}%
holds.

If 
\begin{equation*}
\ell ^{2}\left( \mathbb{K}\right) :=\left\{ \mathbf{a}=\left( a_{i}\right)
_{i\in \mathbb{N}}\subset \mathbb{K}\left\vert \sum_{i=1}^{\infty
}\left\vert a_{i}\right\vert ^{2}\right. <\infty \right\} ,
\end{equation*}%
where $\mathbb{K}=\mathbb{C}$ or $\mathbb{K}=\mathbb{R}$, is the Hilbert
space of all real or complex sequences that are $2-$summable and $\mathbf{%
\lambda }=\left( \lambda _{i}\right) _{i\in \mathbb{N}}\in \ell ^{2}\left( 
\mathbb{K}\right) ,$ then the series $\sum_{i=1}^{\infty }\lambda _{i}e_{i}$
is convergent in $H$ and if $y:=\sum_{i=1}^{\infty }\lambda _{i}e_{i}\in H,$
then $\left\Vert y\right\Vert =\left( \sum_{i=1}^{\infty }\left\vert \lambda
_{i}\right\vert ^{2}\right) ^{\frac{1}{2}}.$

We may state the following result.

\begin{theorem}
\label{t5.1}Let $\left( H;\left\langle \cdot ,\cdot \right\rangle \right) $
be an infinite dimensional Hilbert space over the real or complex number
field $\mathbb{K}$, $\left( e_{i}\right) _{i\in \mathbb{N}}$ is an
orthornormal family in $H,$ $\mathbf{\lambda }=\left( \lambda _{i}\right)
_{i\in \mathbb{N}}\in \ell ^{2}\left( \mathbb{K}\right) ,$ $\mathbf{\lambda }%
\neq 0$ and $r>0.$ If $x\in H$ is such that%
\begin{equation}
\left\Vert x-\sum_{i=1}^{\infty }\lambda _{i}e_{i}\right\Vert \leq r,
\label{5.2}
\end{equation}%
then we have the inequality%
\begin{equation}
0\leq \left\Vert x\right\Vert -\left( \sum_{i=1}^{\infty }\left\vert
\left\langle x,e_{i}\right\rangle \right\vert ^{2}\right) ^{\frac{1}{2}}\leq 
\frac{1}{2}\cdot \frac{r^{2}}{\left( \sum_{i=1}^{\infty }\left\vert \lambda
_{i}\right\vert ^{2}\right) ^{\frac{1}{2}}}.  \label{5.3}
\end{equation}%
The constant $\frac{1}{2}$ is best possible in (\ref{5.3}) in the sense that
it cannot be replaced by a smaller constant.
\end{theorem}

\begin{proof}
Let $a:=\sum_{i=1}^{\infty }\lambda _{i}e_{i}\in H.$ Then by Theorem \ref%
{t2.1}, we have 
\begin{equation*}
\left\Vert x\right\Vert \left\Vert \sum_{i=1}^{\infty }\lambda
_{i}e_{i}\right\Vert -\left\vert \sum_{i=1}^{\infty }\bar{\lambda}%
_{i}\left\langle x,e_{i}\right\rangle \right\vert \leq \frac{1}{2}r^{2},
\end{equation*}%
giving%
\begin{equation}
\left\Vert x\right\Vert \left( \sum_{i=1}^{\infty }\left\vert \lambda
_{i}\right\vert ^{2}\right) ^{\frac{1}{2}}\leq \frac{1}{2}r^{2}+\left\vert
\sum_{i=1}^{\infty }\bar{\lambda}_{i}\left\langle x,e_{i}\right\rangle
\right\vert ,  \label{5.4}
\end{equation}%
since%
\begin{equation*}
\left\Vert \sum_{i=1}^{\infty }\lambda _{i}e_{i}\right\Vert =\left(
\sum_{i=1}^{\infty }\left\vert \lambda _{i}\right\vert ^{2}\right) ^{\frac{1%
}{2}}.
\end{equation*}%
Using the Cauchy-Bunyakovsky-Schwarz inequality, we may state that%
\begin{equation}
\left\vert \sum_{i=1}^{\infty }\bar{\lambda}_{i}\left\langle
x,e_{i}\right\rangle \right\vert \leq \left( \sum_{i=1}^{\infty }\left\vert
\lambda _{i}\right\vert ^{2}\right) ^{\frac{1}{2}}\left( \sum_{i=1}^{\infty
}\left\vert \left\langle x,e_{i}\right\rangle \right\vert ^{2}\right) ^{%
\frac{1}{2}},  \label{5.5}
\end{equation}%
and thus, by (\ref{5.4}) and (\ref{5.5}), we may state that%
\begin{equation*}
\left\Vert x\right\Vert \left( \sum_{i=1}^{\infty }\left\vert \lambda
_{i}\right\vert ^{2}\right) ^{\frac{1}{2}}\leq \frac{1}{2}r^{2}+\left(
\sum_{i=1}^{\infty }\left\vert \lambda _{i}\right\vert ^{2}\right) ^{\frac{1%
}{2}}\left( \sum_{i=1}^{\infty }\left\vert \left\langle x,e_{i}\right\rangle
\right\vert ^{2}\right) ^{\frac{1}{2}},
\end{equation*}%
from where we get the desired inequality in (\ref{5.3}).

The best constant, follows by Theorem \ref{t2.1} on choosing $\left(
e_{i}\right) _{i\in \mathbb{N}}=\left\{ e\right\} ,$ with $\left\Vert
e\right\Vert =1$ and we omit the details.
\end{proof}

\begin{remark}
Under the assumptions of Theorem \ref{t5.1}, and if we multiply by $%
\left\Vert x\right\Vert +\left( \sum_{i=1}^{\infty }\left\vert \left\langle
x,e_{i}\right\rangle \right\vert ^{2}\right) ^{\frac{1}{2}}>0,$ we deduce
from (\ref{5.3}), that%
\begin{align}
0& \leq \left\Vert x\right\Vert ^{2}-\sum_{i=1}^{\infty }\left\vert
\left\langle x,e_{i}\right\rangle \right\vert ^{2}  \label{5.6} \\
& \leq \frac{1}{2}\cdot \frac{r^{2}\left( \left\Vert x\right\Vert +\left(
\sum_{i=1}^{\infty }\left\vert \left\langle x,e_{i}\right\rangle \right\vert
^{2}\right) ^{\frac{1}{2}}\right) }{\left( \sum_{i=1}^{\infty }\left\vert
\lambda _{i}\right\vert ^{2}\right) ^{\frac{1}{2}}}  \notag \\
& \leq \frac{r^{2}\left\Vert x\right\Vert }{\left( \sum_{i=1}^{\infty
}\left\vert \lambda _{i}\right\vert ^{2}\right) ^{\frac{1}{2}}},  \notag
\end{align}%
where for the last inequality, we have used Bessel's inequality%
\begin{equation*}
\left( \sum_{i=1}^{\infty }\left\vert \left\langle x,e_{i}\right\rangle
\right\vert ^{2}\right) ^{\frac{1}{2}}\leq \left\Vert x\right\Vert ,\ \ \ \
x\in H.
\end{equation*}
\end{remark}

The following result also holds.

\begin{theorem}
\label{t5.2}Assume that $\left( H;\left\langle \cdot ,\cdot \right\rangle
\right) $ and $\left( e_{i}\right) _{i\in \mathbb{N}}$ are as in Theorem \ref%
{t5.1}. If $\mathbf{\Gamma }=\left( \Gamma _{i}\right) _{i\in \mathbb{N}},$ $%
\mathbf{\gamma }=\left( \gamma _{i}\right) _{i\in \mathbb{N}}\in \ell
^{2}\left( \mathbb{K}\right) ,$ with $\mathbf{\Gamma }\neq \pm \mathbf{%
\gamma }$, and $x\in H$ with the property that, either%
\begin{equation}
\left\Vert x-\sum_{i=1}^{\infty }\frac{\Gamma _{i}+\gamma _{i}}{2}\cdot
e_{i}\right\Vert \leq \frac{1}{2}\left( \sum_{i=1}^{\infty }\left\vert
\Gamma _{i}-\gamma _{i}\right\vert ^{2}\right) ^{\frac{1}{2}}  \label{5.7}
\end{equation}%
or, equivalently,%
\begin{equation}
\func{Re}\left\langle \sum_{i=1}^{\infty }\Gamma
_{i}e_{i}-x,x-\sum_{i=1}^{\infty }\gamma _{i}e_{i}\right\rangle \geq 0
\label{5.8}
\end{equation}%
holds, then we have the inequality%
\begin{equation}
0\leq \left\Vert x\right\Vert -\left( \sum_{i=1}^{\infty }\left\vert
\left\langle x,e_{i}\right\rangle \right\vert ^{2}\right) ^{\frac{1}{2}}\leq 
\frac{1}{4}\cdot \frac{\sum_{i=1}^{\infty }\left\vert \Gamma _{i}-\gamma
_{i}\right\vert ^{2}}{\left( \sum_{i=1}^{\infty }\left\vert \Gamma
_{i}+\gamma _{i}\right\vert ^{2}\right) ^{\frac{1}{2}}}.  \label{5.9}
\end{equation}%
The constant $\frac{1}{4}$ is best possible in the sense that it cannot be
replaced by a smaller constant.
\end{theorem}

\begin{proof}
Since $\mathbf{\Gamma },$ $\mathbf{\gamma }\in \ell ^{2}\left( \mathbb{K}%
\right) ,$ then we have that $\frac{1}{2}\left( \mathbf{\Gamma }\pm \mathbf{%
\gamma }\right) \in \ell ^{2}\left( \mathbb{K}\right) ,$ showing that the
series%
\begin{equation*}
\sum_{i=1}^{\infty }\left\vert \frac{\Gamma _{i}+\gamma _{i}}{2}\right\vert
^{2},\ \ \ \sum_{i=1}^{\infty }\left\vert \frac{\Gamma _{i}-\gamma _{i}}{2}%
\right\vert ^{2}
\end{equation*}%
are convergent. In addition, the series $\sum_{i=1}^{\infty }\Gamma
_{i}e_{i} $, $\sum_{i=1}^{\infty }\gamma _{i}e_{i}$ and $\sum_{i=1}^{\infty }%
\frac{\Gamma _{i}+\gamma _{i}}{2}e_{i}$ are also convergent in the Hilbert
space $H.$

The equivalence of the conditions (\ref{5.7}) and (\ref{5.8}) follows by the
fact that, in an inner product space we have, for $x,z,Z\in H,$ $\func{Re}%
\left\langle Z-x,x-z\right\rangle \geq 0$ is equivalent to $\left\Vert x-%
\frac{z+Z}{2}\right\Vert \leq \frac{1}{2}\left\Vert Z-z\right\Vert ,$ and we
omit the details.

Now, we observe that the inequality (\ref{5.9}) follows from Theorem \ref%
{t5.1} on choosing $\lambda _{i}=\frac{\Gamma _{i}+\gamma _{i}}{2},$ $i\in 
\mathbb{N}$ and $r=\frac{1}{2}\left( \sum_{i=1}^{\infty }\left\vert \Gamma
_{i}-\gamma _{i}\right\vert ^{2}\right) ^{\frac{1}{2}}.$

The fact that $\frac{1}{4}$ is the best possible constant in (\ref{5.9})
follows from Theorem \ref{t2.2}, and we omit the details.
\end{proof}

\begin{remark}
With the assumptions of Theorem \ref{t5.2}, we have 
\begin{align}
0& \leq \left\Vert x\right\Vert ^{2}-\sum_{i=1}^{\infty }\left\vert
\left\langle x,e_{i}\right\rangle \right\vert ^{2}  \label{5.10} \\
& \leq \frac{1}{4}\cdot \frac{\sum_{i=1}^{\infty }\left\vert \Gamma
_{i}-\gamma _{i}\right\vert ^{2}}{\left( \sum_{i=1}^{\infty }\left\vert
\Gamma _{i}+\gamma _{i}\right\vert ^{2}\right) ^{\frac{1}{2}}}\left[
\left\Vert x\right\Vert +\left( \sum_{i=1}^{\infty }\left\vert \left\langle
x,e_{i}\right\rangle \right\vert ^{2}\right) ^{\frac{1}{2}}\right]  \notag \\
& \leq \frac{1}{2}\cdot \frac{\sum_{i=1}^{\infty }\left\vert \Gamma
_{i}-\gamma _{i}\right\vert ^{2}}{\left( \sum_{i=1}^{\infty }\left\vert
\Gamma _{i}+\gamma _{i}\right\vert ^{2}\right) ^{\frac{1}{2}}}\left\Vert
x\right\Vert .  \notag
\end{align}
\end{remark}

For some recent results related to Bessel inequality, see \cite{HXC}, \cite%
{SSD01}, \cite{SSDJS}, and \cite{GH1}.

\section{Some Gr\"{u}ss Type Inequalities for Orthonormal Families}

The following result holds.

\begin{theorem}
\label{t6.1}Let $\left( H;\left\langle \cdot ,\cdot \right\rangle \right) $
be an infinite dimensional Hilbert space over the real or complex number
field $\mathbb{K}$ and $\left( e_{i}\right) _{i\in \mathbb{N}}$ an
orthornormal family in $H.$ If $\mathbf{\lambda }=\left( \lambda _{i}\right)
_{i\in \mathbb{N}},\ \mathbf{\mu }=\left( \mu _{i}\right) _{i\in \mathbb{N}%
}\in \ell ^{2}\left( \mathbb{K}\right) ,$ $\mathbf{\lambda }$, $\mathbf{\mu }%
\neq 0,$ $r_{1},r_{2}>0$ and $x,y\in H$ are such that%
\begin{equation}
\left\Vert x-\sum_{i=1}^{\infty }\lambda _{i}e_{i}\right\Vert \leq r_{1},\ \
\ \left\Vert y-\sum_{i=1}^{\infty }\mu _{i}e_{i}\right\Vert \leq r_{2},
\label{6.1}
\end{equation}%
then we have the inequality%
\begin{multline}
\left\vert \left\langle x,y\right\rangle -\sum_{i=1}^{\infty }\left\langle
x,e_{i}\right\rangle \left\langle e_{i},y\right\rangle \right\vert
\label{6.2} \\
\leq \frac{1}{2}r_{1}r_{2}\frac{\left[ \left\Vert x\right\Vert +\left(
\sum_{i=1}^{\infty }\left\vert \left\langle x,e_{i}\right\rangle \right\vert
^{2}\right) ^{\frac{1}{2}}\right] ^{\frac{1}{2}}\left[ \left\Vert
y\right\Vert +\left( \sum_{i=1}^{\infty }\left\vert \left\langle
y,e_{i}\right\rangle \right\vert ^{2}\right) ^{\frac{1}{2}}\right] ^{\frac{1%
}{2}}}{\left( \sum_{i=1}^{\infty }\left\vert \lambda _{i}\right\vert
^{2}\right) ^{\frac{1}{4}}\left( \sum_{i=1}^{\infty }\left\vert \mu
_{i}\right\vert ^{2}\right) ^{\frac{1}{4}}}
\end{multline}%
\begin{equation*}
\leq r_{1}r_{2}\frac{\left\Vert x\right\Vert ^{\frac{1}{2}}\left\Vert
y\right\Vert ^{\frac{1}{2}}}{\left( \sum_{i=1}^{\infty }\left\vert \lambda
_{i}\right\vert ^{2}\right) ^{\frac{1}{4}}\left( \sum_{i=1}^{\infty
}\left\vert \mu _{i}\right\vert ^{2}\right) ^{\frac{1}{4}}}.
\end{equation*}
\end{theorem}

\begin{proof}
Apply Schwarz's inequality for the vectors $x-\sum_{i=1}^{\infty
}\left\langle x,e_{i}\right\rangle e_{i},$ $y-\sum_{i=1}^{\infty
}\left\langle y,e_{i}\right\rangle e_{i},$ to get%
\begin{multline}
\left\vert \left\langle x-\sum_{i=1}^{\infty }\left\langle
x,e_{i}\right\rangle e_{i},y-\sum_{i=1}^{\infty }\left\langle
y,e_{i}\right\rangle e_{i}\right\rangle \right\vert ^{2}  \label{6.3} \\
\leq \left\Vert x-\sum_{i=1}^{\infty }\left\langle x,e_{i}\right\rangle
e_{i}\right\Vert ^{2}\left\Vert y-\sum_{i=1}^{\infty }\left\langle
y,e_{i}\right\rangle e_{i}\right\Vert ^{2}.
\end{multline}%
Since%
\begin{equation*}
\left\langle x-\sum_{i=1}^{\infty }\left\langle x,e_{i}\right\rangle
e_{i},y-\sum_{i=1}^{\infty }\left\langle y,e_{i}\right\rangle
e_{i}\right\rangle =\left\langle x,y\right\rangle -\sum_{i=1}^{\infty
}\left\langle x,e_{i}\right\rangle \left\langle e_{i},y\right\rangle
\end{equation*}%
and%
\begin{equation*}
\left\Vert x-\sum_{i=1}^{\infty }\left\langle x,e_{i}\right\rangle
e_{i}\right\Vert ^{2}=\left\Vert x\right\Vert ^{2}-\sum_{i=1}^{\infty
}\left\vert \left\langle x,e_{i}\right\rangle \right\vert ^{2},
\end{equation*}%
then, by (\ref{6.3}) and (\ref{5.6}) applied for $x$ and $y,$ we deduce the
desired inequality (\ref{6.2}).
\end{proof}

Finally we may state the following theorem.

\begin{theorem}
\label{t6.2}Assume that $\left( H;\left\langle \cdot ,\cdot \right\rangle
\right) $ and $\left( e_{i}\right) _{i\in \mathbb{N}}$ are as in Theorem \ref%
{t6.1}. If $\mathbf{\Gamma }=\left( \Gamma _{i}\right) _{i\in \mathbb{N}},$ $%
\mathbf{\Gamma }=\left( \Gamma _{i}\right) _{i\in \mathbb{N}},\phi =\left( 
\mathbf{\phi }_{i}\right) _{i\in \mathbb{N}},\mathbf{\Phi }=\left( \Phi
_{i}\right) _{i\in \mathbb{N}}\in \ell ^{2}\left( \mathbb{K}\right) ,$ with $%
\mathbf{\Gamma }\neq \pm \mathbf{\gamma }$, $\mathbf{\Phi }\neq \pm \mathbf{%
\phi }$, and $x,y\in H$ are such that, either%
\begin{eqnarray}
\left\Vert x-\sum_{i=1}^{\infty }\frac{\Gamma _{i}+\gamma _{i}}{2}\cdot
e_{i}\right\Vert &\leq &\frac{1}{2}\left( \sum_{i=1}^{\infty }\left\vert
\Gamma _{i}-\gamma _{i}\right\vert ^{2}\right) ^{\frac{1}{2}},  \label{6.4}
\\
\left\Vert y-\sum_{i=1}^{\infty }\frac{\Phi _{i}+\phi _{i}}{2}\cdot
e_{i}\right\Vert &\leq &\frac{1}{2}\left( \sum_{i=1}^{\infty }\left\vert
\Phi _{i}-\phi _{i}\right\vert ^{2}\right) ^{\frac{1}{2}},  \notag
\end{eqnarray}%
or, equivalently,%
\begin{eqnarray}
\func{Re}\left\langle \sum_{i=1}^{\infty }\Gamma
_{i}e_{i}-x,x-\sum_{i=1}^{\infty }\gamma _{i}e_{i}\right\rangle &\geq &0,
\label{6.5} \\
\func{Re}\left\langle \sum_{i=1}^{\infty }\Phi
_{i}e_{i}-y,y-\sum_{i=1}^{\infty }\phi _{i}e_{i}\right\rangle &\geq &0,
\end{eqnarray}%
holds, then we have the inequality%
\begin{align}
& \left\vert \left\langle x,y\right\rangle -\sum_{i=1}^{\infty }\left\langle
x,e_{i}\right\rangle \left\langle e_{i},y\right\rangle \right\vert
\label{6.6} \\
& \leq \frac{1}{4}\cdot \left( \sum_{i=1}^{\infty }\left\vert \Phi _{i}-\phi
_{i}\right\vert ^{2}\right) ^{\frac{1}{2}}\left( \sum_{i=1}^{\infty
}\left\vert \Gamma _{i}-\gamma _{i}\right\vert ^{2}\right) ^{\frac{1}{2}} 
\notag \\
& \qquad \qquad \times \frac{\left[ \left\Vert x\right\Vert +\left(
\sum_{i=1}^{\infty }\left\vert \left\langle x,e_{i}\right\rangle \right\vert
^{2}\right) ^{\frac{1}{2}}\right] ^{\frac{1}{2}}\left[ \left\Vert
y\right\Vert +\left( \sum_{i=1}^{\infty }\left\vert \left\langle
y,e_{i}\right\rangle \right\vert ^{2}\right) ^{\frac{1}{2}}\right] ^{\frac{1%
}{2}}}{\left( \sum_{i=1}^{\infty }\left\vert \Phi _{i}+\phi _{i}\right\vert
^{2}\right) ^{\frac{1}{4}}\left( \sum_{i=1}^{\infty }\left\vert \Gamma
_{i}+\gamma _{i}\right\vert ^{2}\right) ^{\frac{1}{4}}}  \notag \\
& \leq \frac{1}{2}\cdot \frac{\left( \sum_{i=1}^{\infty }\left\vert \Phi
_{i}-\phi _{i}\right\vert ^{2}\right) ^{\frac{1}{2}}\left(
\sum_{i=1}^{\infty }\left\vert \Gamma _{i}-\gamma _{i}\right\vert
^{2}\right) ^{\frac{1}{2}}}{\left( \sum_{i=1}^{\infty }\left\vert \Phi
_{i}+\phi _{i}\right\vert ^{2}\right) ^{\frac{1}{4}}\left(
\sum_{i=1}^{\infty }\left\vert \Gamma _{i}+\gamma _{i}\right\vert
^{2}\right) ^{\frac{1}{4}}}\left\Vert x\right\Vert ^{\frac{1}{2}}\left\Vert
y\right\Vert ^{\frac{1}{2}}.  \notag
\end{align}
\end{theorem}

The proof follows by (\ref{6.3}) and by (\ref{5.10}) applied for $x$ and $y.$
We omit the details.

\section{Integral Inequalities}

Let $\left( \Omega ,\Sigma ,\mu \right) $ be a measure space consisting of a
set $\Omega ,$ a $\sigma -$algebra of parts $\Sigma $ and a countably
additive and positive measure $\mu $ on $\Sigma $ with values in $\mathbb{R}%
\cup \left\{ \infty \right\} .$ Let $\rho \geq 0$ be a $\mu -$measurable
function on $\Omega $ with $\int_{\Omega }\rho \left( s\right) d\mu \left(
s\right) =1.$ Denote by $L_{\rho }^{2}\left( \Omega ,\mathbb{K}\right) $ the
Hilbert space of all real or complex valued functions defined on $\Omega $
and $2-\rho -$integrable on $\Omega ,$ i.e.,%
\begin{equation}
\int_{\Omega }\rho \left( s\right) \left\vert f\left( s\right) \right\vert
^{2}d\mu \left( s\right) <\infty .  \label{7.1}
\end{equation}%
It is obvious that the following inner product%
\begin{equation}
\left\langle f,g\right\rangle _{\rho }:=\int_{\Omega }\rho \left( s\right)
f\left( s\right) \overline{g\left( s\right) }d\mu \left( s\right) ,
\label{7.2}
\end{equation}%
generates the norm 
\begin{equation*}
\left\Vert f\right\Vert _{\rho }:=\left( \int_{\Omega }\rho \left( s\right)
\left\vert f\left( s\right) \right\vert ^{2}d\mu \left( s\right) \right) ^{%
\frac{1}{2}}
\end{equation*}%
of $L_{\rho }^{2}\left( \Omega ,\mathbb{K}\right) ,$ and all the above
results may be stated for integrals.

It is important to observe that, if%
\begin{equation}
\func{Re}\left[ f\left( s\right) \overline{g\left( s\right) }\right] \geq 0%
\text{ \ \ for \ }\mu -\text{a.e. \ }s\in \Omega ,  \label{7.3}
\end{equation}%
then, obviously,%
\begin{align}
\func{Re}\left\langle f,g\right\rangle _{\rho }& =\func{Re}\left[
\int_{\Omega }\rho \left( s\right) f\left( s\right) \overline{g\left(
s\right) }d\mu \left( s\right) \right]  \label{7.4} \\
& =\int_{\Omega }\rho \left( s\right) \func{Re}\left[ f\left( s\right) 
\overline{g\left( s\right) }\right] d\mu \left( s\right) \geq 0.  \notag
\end{align}%
The reverse is evidently not true in general.

Moreover, if the space is real, i.e., $\mathbb{K}=\mathbb{R}$, then a
sufficient condition for (\ref{7.4}) to hold is:%
\begin{equation}
f\left( s\right) \geq 0,\ \ g\left( s\right) \geq 0\text{ \ \ for \ }\mu -%
\text{a.e. \ }s\in \Omega .  \label{7.5}
\end{equation}

We provide now, by the use of certain results obtained in Section \ref{s2},
some integral inequalities that may be used in practical applications.

\begin{proposition}
\label{p7.1}Let $f,g\in L_{\rho }^{2}\left( \Omega ,\mathbb{K}\right) $ and $%
r>0$ with the property that%
\begin{equation}
\left\vert f\left( s\right) -g\left( s\right) \right\vert \leq r\text{ \ \
for \ }\mu -\text{a.e. \ }s\in \Omega .  \label{7.6}
\end{equation}%
Then we have the inequalities%
\begin{align}
0& \leq \left[ \int_{\Omega }\rho \left( s\right) \left\vert f\left(
s\right) \right\vert ^{2}d\mu \left( s\right) \int_{\Omega }\rho \left(
s\right) \left\vert g\left( s\right) \right\vert ^{2}d\mu \left( s\right) %
\right] ^{\frac{1}{2}}  \label{7.7} \\
& \qquad \qquad -\left\vert \int_{\Omega }\rho \left( s\right) f\left(
s\right) \overline{g\left( s\right) }d\mu \left( s\right) \right\vert  \notag
\\
& \leq \left[ \int_{\Omega }\rho \left( s\right) \left\vert f\left( s\right)
\right\vert ^{2}d\mu \left( s\right) \int_{\Omega }\rho \left( s\right)
\left\vert g\left( s\right) \right\vert ^{2}d\mu \left( s\right) \right] ^{%
\frac{1}{2}}  \notag \\
& \qquad \qquad -\left\vert \int_{\Omega }\rho \left( s\right) \func{Re}%
\left[ f\left( s\right) \overline{g\left( s\right) }\right] d\mu \left(
s\right) \right\vert  \notag \\
& \leq \left[ \int_{\Omega }\rho \left( s\right) \left\vert f\left( s\right)
\right\vert ^{2}d\mu \left( s\right) \int_{\Omega }\rho \left( s\right)
\left\vert g\left( s\right) \right\vert ^{2}d\mu \left( s\right) \right] ^{%
\frac{1}{2}}  \notag \\
& \qquad \qquad -\int_{\Omega }\rho \left( s\right) \func{Re}\left[ f\left(
s\right) \overline{g\left( s\right) }\right] d\mu \left( s\right)  \notag \\
& \leq \frac{1}{2}r^{2}.  \notag
\end{align}%
The constant $\frac{1}{2}$ is best possible in (\ref{7.7}).
\end{proposition}

The proof follows by Theorem \ref{t2.1}, and we omit the details.

\begin{proposition}
\label{p7.2}Let $f,g\in L_{\rho }^{2}\left( \Omega ,\mathbb{K}\right) $ and $%
\gamma ,\Gamma \in \mathbb{K}$ so that $\Gamma \neq -\gamma ,\gamma $ and%
\begin{equation}
\func{Re}\left[ \left( \Gamma g\left( s\right) -f\left( s\right) \right)
\left( \overline{f\left( s\right) }-\overline{\gamma }\overline{g\left(
s\right) }\right) \right] \geq 0\text{ \ \ for \ }\mu -\text{a.e. \ }s\in
\Omega .  \label{7.8}
\end{equation}%
Then we have the inequalities%
\begin{align}
0& \leq \left[ \int_{\Omega }\rho \left( s\right) \left\vert f\left(
s\right) \right\vert ^{2}d\mu \left( s\right) \int_{\Omega }\rho \left(
s\right) \left\vert g\left( s\right) \right\vert ^{2}d\mu \left( s\right) %
\right] ^{\frac{1}{2}}  \label{7.9} \\
& \qquad \qquad -\left\vert \int_{\Omega }\rho \left( s\right) f\left(
s\right) \overline{g\left( s\right) }d\mu \left( s\right) \right\vert  \notag
\\
& \leq \left[ \int_{\Omega }\rho \left( s\right) \left\vert f\left( s\right)
\right\vert ^{2}d\mu \left( s\right) \int_{\Omega }\rho \left( s\right)
\left\vert g\left( s\right) \right\vert ^{2}d\mu \left( s\right) \right] ^{%
\frac{1}{2}}  \notag \\
& \qquad \qquad -\left\vert \func{Re}\left[ \frac{\bar{\Gamma}+\bar{\gamma}}{%
\left\vert \Gamma +\gamma \right\vert }\int_{\Omega }\rho \left( s\right)
f\left( s\right) \overline{g\left( s\right) }d\mu \left( s\right) \right]
\right\vert  \notag \\
& \leq \left[ \int_{\Omega }\rho \left( s\right) \left\vert f\left( s\right)
\right\vert ^{2}d\mu \left( s\right) \int_{\Omega }\rho \left( s\right)
\left\vert g\left( s\right) \right\vert ^{2}d\mu \left( s\right) \right] ^{%
\frac{1}{2}}  \notag \\
& \qquad \qquad -\func{Re}\left[ \frac{\bar{\Gamma}+\bar{\gamma}}{\left\vert
\Gamma +\gamma \right\vert }\int_{\Omega }\rho \left( s\right) f\left(
s\right) \overline{g\left( s\right) }d\mu \left( s\right) \right]  \notag \\
& \leq \frac{1}{4}\cdot \frac{\left\vert \Gamma -\gamma \right\vert ^{2}}{%
\left\vert \Gamma +\gamma \right\vert }\int_{\Omega }\rho \left( s\right)
\left\vert g\left( s\right) \right\vert ^{2}d\mu \left( s\right) .  \notag
\end{align}%
The constant $\frac{1}{4}$ is best possible.
\end{proposition}

\begin{remark}
If the space is real and we assume, for $M>m>0,$ that%
\begin{equation}
mg\left( s\right) \leq f\left( s\right) \leq Mg\left( s\right) \text{ \ \
for \ }\mu -\text{a.e. \ }s\in \Omega ,  \label{7.10}
\end{equation}%
then by (\ref{7.9}) we deduce the inequality:%
\begin{eqnarray}
0 &\leq &\left[ \int_{\Omega }\rho \left( s\right) \left\vert f\left(
s\right) \right\vert ^{2}d\mu \left( s\right) \int_{\Omega }\rho \left(
s\right) \left\vert g\left( s\right) \right\vert ^{2}d\mu \left( s\right) %
\right] ^{\frac{1}{2}}  \label{7.11} \\
&&\qquad \qquad -\left\vert \int_{\Omega }\rho \left( s\right) f\left(
s\right) \overline{g\left( s\right) }d\mu \left( s\right) \right\vert  \notag
\\
&\leq &\frac{1}{4}\cdot \frac{\left( M-m\right) ^{2}}{M+m}\int_{\Omega }\rho
\left( s\right) \left\vert g\left( s\right) \right\vert ^{2}d\mu \left(
s\right) .  \notag
\end{eqnarray}%
The constant $\frac{1}{4}$ is best possible.
\end{remark}

The following reverse of the triangle inequality for integrals holds.

\begin{proposition}
\label{p.7.3}Assume that the functions $f,g\in L_{\rho }^{2}\left( \Omega ,%
\mathbb{K}\right) $ satisfy (\ref{7.10}). Then we have the inequality%
\begin{eqnarray}
0 &\leq &\left( \int_{\Omega }\rho \left( s\right) \left\vert f\left(
s\right) \right\vert ^{2}d\mu \left( s\right) \right) ^{1/2}+\left(
\int_{\Omega }\rho \left( s\right) \left\vert g\left( s\right) \right\vert
^{2}d\mu \left( s\right) \right) ^{1/2}  \label{7.12} \\
&&-\left( \int_{\Omega }\rho \left( s\right) \left\vert f\left( s\right)
+g\left( s\right) \right\vert ^{2}d\mu \left( s\right) \right) ^{1/2}  \notag
\\
&\leq &\frac{\sqrt{2}}{2}\cdot \frac{\left( M-m\right) }{\sqrt{M+m}}\left(
\int_{\Omega }\rho \left( s\right) \left\vert g\left( s\right) \right\vert
^{2}d\mu \left( s\right) \right) ^{1/2}.  \notag
\end{eqnarray}
\end{proposition}

The proof follows by Proposition \ref{p2.4}.

By making use of Theorem \ref{t4.3}, we may also state

\begin{proposition}
\label{p7.3}Let $f,g,h\in L_{\rho }^{2}\left( \Omega ,\mathbb{K}\right) $ be
so that $\int_{\Omega }\rho \left( s\right) \left\vert h\left( s\right)
\right\vert ^{2}d\mu \left( s\right) =1.$ Suppose also that $a,A,b,B\in 
\mathbb{K}$ with $A\neq \pm a,B\neq \pm b$ and%
\begin{eqnarray*}
\func{Re}\left[ \left( Ah\left( s\right) -f\left( s\right) \right) \left( 
\overline{f\left( s\right) }-\overline{a}\overline{h\left( s\right) }\right) %
\right] &\geq &0\text{,} \\
\func{Re}\left[ \left( Bh\left( s\right) -g\left( s\right) \right) \left( 
\overline{g\left( s\right) }-\overline{b}\overline{h\left( s\right) }\right) %
\right] &\geq &0\text{ \ \ for \ }\mu -\text{a.e. \ }s\in \Omega .
\end{eqnarray*}%
Then we have the inequality%
\begin{eqnarray*}
&&\left\vert \int_{\Omega }\rho \left( s\right) f\left( s\right) \overline{%
g\left( s\right) }d\mu \left( s\right) -\int_{\Omega }\rho \left( s\right)
f\left( s\right) \overline{h\left( s\right) }d\mu \left( s\right)
\int_{\Omega }\rho \left( s\right) h\left( s\right) \overline{g\left(
s\right) }d\mu \left( s\right) \right\vert \\
&\leq &\frac{1}{4}\cdot \frac{\left\vert A-a\right\vert \left\vert
B-b\right\vert }{\sqrt{\left\vert A+a\right\vert \left\vert B+b\right\vert }}
\\
&&\times \sqrt{\left( \int_{\Omega }\rho \left( s\right) \left\vert f\left(
s\right) \right\vert ^{2}d\mu \left( s\right) \right) ^{1/2}+\left\vert
\int_{\Omega }\rho \left( s\right) f\left( s\right) \overline{h\left(
s\right) }d\mu \left( s\right) \right\vert } \\
&&\times \sqrt{\left( \int_{\Omega }\rho \left( s\right) \left\vert g\left(
s\right) \right\vert ^{2}d\mu \left( s\right) \right) ^{1/2}+\left\vert
\int_{\Omega }\rho \left( s\right) g\left( s\right) \overline{h\left(
s\right) }d\mu \left( s\right) \right\vert }.
\end{eqnarray*}%
The constant $\frac{1}{4}$ is best possible.
\end{proposition}

\begin{remark}
All the other inequalities in Sections \ref{s3} -- \ref{s6} may be used in a
similar manner to obtain the corresponding integral inequalities. We omit
the details.
\end{remark}


\begin{thebibliography}{99}
\bibitem{HXC} X.H. CAO, Bessel sequences in a Hilbert space. \textit{%
Gongcheng Shuxue Xuebao} 1\textbf{7} (2000), no. 2, 92--98.

\bibitem{ADR} A. De ROSSI, A strengthened Cauchy-Schwarz inequality for
biorthogonal wavelets. \textit{Math. Inequal. Appl.} \textbf{2} (1999), no.
2, 263--282.

\bibitem{JBDFTM} J. B. DIAZ and F. T. METCALF, A complementary triangle
inequality in Hilbert and Banach spaces. \textit{Proc. Amer. Math. Soc}. 
\textbf{17} (1966), 88--97.

\bibitem{SSD0} S.S. DRAGOMIR, A generalization of Gr\"{u}ss inequality in
inner product spaces and applications. \textit{J. Math. Anal. Appl.} \textbf{%
237} (1999), no. 1, 74--82.

\bibitem{SSD01} S.S. DRAGOMIR, A note on Bessel's inequality, \textit{%
Austral. Math. Soc. Gaz.} \textbf{28} (2001), no. 5, 246--248.

\bibitem{SSD00} S.S. DRAGOMIR, Some Gr\"{u}ss type inequalities in inner
product spaces, \textit{J. Inequal. Pure \& Appl. Math.}, \textbf{4}(2003),
No. 2, Article 42, [\texttt{On line: http://jipam.vu.edu.au/v4n2/032\_03.html%
}]

\bibitem{SSD1} S.S. DRAGOMIR, A counterpart of Schwarz's inequality in inner
product spaces, \textit{RGMIA Res. Rep. Coll.,} \textbf{6}(2003), \textit{%
Supplement}, Article 18, \texttt{[On line http://rgmia.vu.edu.au/v6(E).html]}

\bibitem{SSD2} S.S. DRAGOMIR, A generalisation of the Cassels and
Grueb-Reinboldt inequalities in inner product spaces, Preprint, \textit{%
Mathematics Ar}$X$\textit{iv}, math.CA/0307130, \texttt{[On line
http://front.math.ucdavis.edu/]}

\bibitem{SSD3} S.S. DRAGOMIR, Some companions of the Gr\"{u}ss inequality in
inner product spaces, \textit{RGMIA Res. Rep. Coll. }\textbf{6}(2003), 
\textit{Supplement}, Article 8, \texttt{[On line
http://rgmia.vu.edu.au/v6(E).html]}

\bibitem{SSD4} S.S. DRAGOMIR, On Bessel and Gr\"{u}ss inequalities for
orthornormal families in inner product spaces, \textit{RGMIA Res. Rep. Coll. 
}\textbf{6}(2003), \textit{Supplement}, Article 12, \texttt{[On line
http://rgmia.vu.edu.au/v6(E).html]}

\bibitem{SSD5} S.S. DRAGOMIR, A counterpart of Bessel's inequality in inner
product spaces and some Gr\"{u}ss type related results, \textit{RGMIA Res.
Rep. Coll. }\textbf{6}(2003), \textit{Supplement}, Article 10, \texttt{[On
line http://rgmia.vu.edu.au/v6(E).html]}

\bibitem{SSD6} S.S. DRAGOMIR, Some new results related to Bessel and Gr\"{u}%
ss inequalities for orthornormal families in inner product spaces, \textit{%
RGMIA Res. Rep. Coll. }\textbf{6}(2003), \textit{Supplement}, Article 13, 
\texttt{[On line http://rgmia.vu.edu.au/v6(E).html]}

\bibitem{SSDa} S.S. DRAGOMIR, Reverses of Schwarz, triangle and Bessel
inequalities in inner product spaces, \textit{Preprint on line: }%
http://www.mathpreprints.com/math/Preprint/Sever/20030828.2/1/?=%
\&coll=Selection

\bibitem{DM} S. S. DRAGOMIR and B. MOND, On the superadditivity and
monotonicity of Schwarz's inequality in inner product spaces. \textit{%
Makedon. Akad. Nauk. Umet. Oddel. Mat.-Tehn. Nauk. Prilozi }\textbf{15}
(1994), no. 2, 5--22 (1996).

\bibitem{SSDJS} S.S. DRAGOMIR and J. S\'{A}NDOR, On Bessel's and Gram's
inequalities in pre-Hilbertian spaces. \textit{Period. Math. Hungar.} 
\textbf{29} (1994), no. 3, 197--205.

\bibitem{GH} H. GUNAWAN, On n-inner products, n-norms, and the
Cauchy-Schwarz inequality. \textit{Sci. Math. Jpn.} \textbf{55} (2002), no.
1, 53--60.

\bibitem{GH1} H. GUNAWAN, A generalization of Bessel's inequality and
Parseval's identity. \textit{Period. Math. Hungar.} \textbf{44} (2002), no.
2, 177--181.

\bibitem{SMK} S.M. KHALEELULLA, On Diaz-Metcalf's complementary triangle
inequality.\textit{\ Kyungpook Math. J.} \textbf{15 }(1975), 9--11.

\bibitem{PMM} P.M. MILI\v{C}I\'{C}, On a complementary inequality of the
triangle inequality (French)\textit{\ Mat. Vesnik}\textbf{\ 41} (1989), no.
2, 83--88.

\bibitem{DKR} D. K. RAO, A triangle inequality for angles in a Hilbert
space. \textit{Rev. Colombiana Mat.} \textbf{10} (1976), no. 3, 95--97.

\bibitem{PFR} P. F. RENAUD, A matrix formulation of Gr\"{u}ss inequality, 
\textit{Linear Algebra Appl.} \textbf{335} (2001), 95--100.

\bibitem{NU} N. UJEVI\'{C}, A generalisation of Gr\"{u}ss inequality in
prehilbertian spaces, \textit{Math. Inequal. \& Appl., }(to appear).
\end{thebibliography}
\end{document}